\documentclass{amsart}
\usepackage{amscd,amsmath}
\usepackage[dvips]{graphicx}
\usepackage{hyperref}
\usepackage[latin1]{inputenc}
\usepackage{array}
\usepackage[all]{xy}

\newcommand{\sing}{{\rm{Sing}}}

\newtheorem{theorem}{\bf Theorem}
\newtheorem{proposition}{\bf Proposition}
\newtheorem{lemma}{\bf Lemma}

\newtheorem{definition}{\bf Definition}

\newtheorem{corollary}{\bf Corollary}
\newtheorem{example}{\bf Example}

\begin{document}

\title[Polar Transformations versus Logarithmic Foliations]{On the degree of Polar Transformations \\  An approach through Logarithmic Foliations}

\author{T. FASSARELLA}
\address{IMPA  \\
Estrada Dona Castorina, 110 --
Jardim Bot\^{a}nico \\
22460-320  RJ Brasil} \email{thiago@impa.br}

\author{J. V. PEREIRA}
\address{IMPA  \\
Estrada Dona Castorina, 110 --
Jardim Bot\^{a}nico \\
22460-320  RJ Brasil} \email{jvp@impa.br}

 \subjclass{}

\begin{abstract}
We investigate the degree of the polar transformations associated
to a certain class of multi-valued homogeneous functions. In particular
we prove that the degree of the pre-image of  generic linear spaces  
by a polar transformation associated to a homogeneous polynomial $F$ is determined
by the zero locus of $F$.  For zero dimensional-dimensional linear spaces this was
conjecture by Dolgachev and proved by Dimca-Papadima using topological arguments. Our methods
are algebro-geometric and rely on the study of the Gauss map of naturally associated logarithmic
foliations. 
\end{abstract}

\maketitle

\section{Introduction}

Given a homogeneous polynomial  $F \in \mathbb C  [x_0,...,x_n]$ one
can naturally associate to it the rational map induced by its linear
system of polars. Explicitly this map can be written as
\begin{eqnarray*}
 \nabla F: \mathbb P^n &\dashrightarrow& \mathbb P^n  \\
  x &\mapsto& \left(\frac{\partial F}{\partial x_0}(x):...:\frac{\partial F}{\partial x_n}(x) \right) \, ,
\end{eqnarray*}
and is the so called  {\bf polar transformation} or {\bf polar map}
of $F$.

The particular case when $\nabla F$ is  a birational map is of
particular interest \cite{ES,EKP,CRS} and in this situation the
polynomial $F$ is said to be {\bf homaloidal}. The classification of
reduced homaloidal polynomials in three variables was carried out by
Dolgachev in \cite{Do00}. It says that $F \in \mathbb
C[x_0,x_1,x_2]$ is a reduced homaloidal polynomial if, and only if,
its  (set theoretical) zero locus $V(F) \subset \mathbb P^2$ has one
of the following forms:
\begin{enumerate}
 \item a smooth conic;
 \item the union of three lines in general position;
 \item the union of a smooth conic and a line tangent to it.
\end{enumerate}

In loc. cit. it is conjectured that the reduceness of $F$ is not
necessary to draw the same conclusion. More precisely it is
conjectured that the degree of $\nabla F$ can be written as a
function of $V(F)$.

  Dimca and  Papadima  \cite{DiPa03}
settled Dolgachev's conjecture by  proving  that for a polynomial $F
\in \mathbb C[x_0,\ldots,x_n]$ the complement $D(F)=\mathbb P^n
\backslash V(F)$ is homotopically equivalent to a  CW-complex
obtained from  $D(F) \cap H$ by attaching $deg(\nabla F)$ cells of
dimension $n$, where $H \subset \mathbb P^n$ is a generic
hyperplane. In particular the degree of $\nabla F$ can be expressed
as
\[
 deg(\nabla F)= {(-1)}^n \chi (D(F) \backslash H) \, .
\]
Their proof is topological and  relies on complex Morse Theory. In
\cite{Dimca}, as well as in \cite{CRS}, the problem of giving an
algebro-geometric proof of Dolgachev's conjecture is raised. Partial
answers have been provided by \cite{KS} and \cite{B}.

The main goal of this  paper is to provided one such
algebro-geometric proof, cf. Theorem \ref{T:main}, by relating the
degree of $\nabla F$ to the degree of the Gauss map of some
naturally associated logarithmic foliations.

Our method allow us also to  deal with the  higher order degrees of
$\nabla F$ --- the degrees of the closure of pre-images of generic
linear subspaces --- and  with more general functions than the
polynomial ones, cf. \S \ref{S:geral}.

 The paper is organized as follows. In \S \ref{S:fol} we recall some basic definitions concerning
holomorphic foliations and their Gauss map and prove Theorem
\ref{T:Gauss} that express the higher order degrees of such Gauss
maps  in terms of the topological degree of the Gauss maps of
generic linear sections of the corresponding  foliations. In \S
\ref{S:log} we study the Gauss maps of logarithmic foliations and
prove that their topological degrees --- under suitable hypotheses
--- can expressed in terms of the top chern class of certain sheaves
of logarithmic differentials. In \S \ref{S:geral} we prove Theorem
\ref{T:polar} that relates the degrees of the polar map with the
ones of a naturally  associated logarithmic foliation. Finally, in
\S \ref{S:main} we prove Theorem \ref{T:main} --- our main result
--- and make a couple of side remarks.
\medskip

\noindent{\bf Acknowledgements:} We first heard about the degree of
polar maps from Francesco Russo. We thank him for calling our
attention to the problem of giving an algebro-geometric  proof of
Dolgachev's conjecture and for his interest in this work. We also
thanks Charles Favre and Dominique Cerveau. The first for enquiring
us about the
 higher order
degrees of polar maps and the second for suggesting that something
like Corollary
 \ref{C:prod} should hold true.

\section{Foliations and Their Gauss Maps}\label{S:fol}

A {\bf codimension one singular holomorphic foliation}, from now on
just  foliation, $\mathcal F$ of a complex manifold $M$ is
determined by a line bundle $\mathcal L$ and an element $\omega \in
\mathrm H^0(M, \Omega^1_M \otimes \mathcal L)$ satisfying
\begin{enumerate}
\item[(i)] $\mathrm{codim} \,\sing(\omega) \ge 2$ where $\sing(\omega)
= \{ x \in M \, \vert \, \omega(x) = 0 \}$;
\item[(ii)] $\omega \wedge d\omega =0$ in $\mathrm H^0(M,\Omega^3_M \otimes \mathcal
L^{\otimes 2}).$
\end{enumerate}

If we drop condition (ii) we obtain the definition of a {\bf
codimension one singular holomorphic distribution} $\mathcal D$  on
$M$. Although we will state  the results of this section for
foliations they  can all be rephrased  for distributions.

The {\bf singular set} of $\mathcal F$, for short $\sing(\mathcal
F)$, is by definition equal to  $\sing(\omega)$. The integrability
condition (ii) determines in an analytic neighborhood of every point
$p \in M \setminus \sing(\mathcal F)$ a holomorphic fibration with
relative tangent sheaf coinciding with the subsheaf of $TM$
determined by the kernel of $\omega$. Analytic continuation of the
fibers of this fibration describes the leaves of $\mathcal F$.

In our study the isolated singularities of $\mathcal F$ will play a
key role. One of the most basic invariants attached to them is their
{\bf multiplicity} $m(\mathcal F,p)$ defined as the intersection
multiplicity at $p$ of the zero section of $\Omega^1_M\otimes
\mathcal L$ with the graph of $\omega$.

In this paper we will focus on the case $M=\mathbb P^n$. The {\bf
degree} of  a foliation of $\mathbb P^n$ is geometrically defined as
the number of tangencies of $\mathcal F$ with a generic line $\ell
\subset \mathbb P^n$. If $\iota: \ell \to \mathbb P^n$ is the
inclusion of such a line then the degree of $\mathcal F$ is the
degree of the zero divisor of the twisted $1$-form  $\iota^*\omega
\in \mathrm H^0(\mathbb \ell, \Omega^1_{\ell} \otimes \mathcal L_{|
\ell})$. Thus the degree of $\mathcal F$ is nothing more than $
\deg(\mathcal L) -2$.

\subsection{The Gauss Map}
The {\bf Gauss map} of a foliation $\mathcal F$ of $\mathbb P^n$ is
the rational map
\begin{eqnarray*}
\mathcal G ( \mathcal F) : \mathbb P^n &\dashrightarrow& \check
{\mathbb P}^n \,  \\
p &\mapsto& T_p \mathcal F
\end{eqnarray*}
where $T_p \mathcal F$ is the projective tangent space of the leaf
of $\mathcal F$ through $p$.

It follows from  Euler's sequence that a $1$-form $\omega \in
\mathrm H^0(\mathbb P^n ,\Omega^1 ( \deg(\mathcal F) + 2) )$ can be
interpreted as  a homogeneous $1$-form on $\mathbb C^{n+1}$, still
denoted by $\omega$,
\[
\omega = \sum_{i=0}^n a_i dx_i
\]
with the coefficients $a_i$ being homogenous polynomials of degree
$\deg(\mathcal F) + 1$ and satisfying Euler's relation $ i_R \omega
= 0,$ where  $i_R$ stands for the interior product with the radial
(or Euler's) vector field $R = \sum_{i=0}^n x_i
\frac{\partial}{\partial x_i}$.

If we interpret $[dx_0: \ldots : dx_n]$ as projective coordinates of
$\check{\mathbb P}^n$ then the Gauss map of the corresponding
 $\mathcal F$  can be written  as   $\mathcal
G(\mathcal F)(p)=[a_0(p): \ldots: a_n(p) ]$.

\subsection{Linear Sections of Foliations}

Assume that $1\le k<n$ and  let  $\iota: \mathbb P^k \to  \mathbb
P^n$ be a linear embedding. If $\iota^* \omega = 0$ then we say that
$\iota(\mathbb P^k)$ is {\bf left invariant} by $\mathcal F$.
Otherwise, after dividing $\iota^* \omega$ (here interpreted as a
$1$-form on $\mathbb C^{k+1}$) by the common factor of its
coefficients, one obtains a foliation $i^* \mathcal F=\mathcal F_{|
\mathbb P^k}$ on $\mathbb P^k$.

Notice that according to our definitions there is only one foliation
of $\mathbb P^1$ and it is  induced by the homogeneous $1$-form $x_0
dx_1- x_1 dx_0$ on $\mathbb C^2$. This odd remark will prove to be
useful when we define the numbers $e_i^k(\mathcal F)$ below. On the
other hand  if $k\ge 2$ and $\iota:\mathbb P^k\to \mathbb P^n$ is
generic then there is no need to divide $\iota^*\omega$: one has
just to apply the following well-known lemma $n-k$ times.

\begin{lemma}\label{L:CLS}
Let $n \ge 3$. If $H \subset \mathbb P^n$ is a generic hyperplane and $\mathcal F$
is a foliation of $\mathbb P^n$  then the degree of $\mathcal F_{|
H}$ is equal to the degree of $\mathcal F$ and, moreover,
\[
\sing(\mathcal F_{|H}) = (\sing(\mathcal F)\cap H ) \cup \mathcal
G(\mathcal F)^{-1}(H)
\]
with $ \mathcal G(\mathcal F)^{-1}(H)$ being finite and all the
corresponding singularities of $\mathcal F_{|H}$ have multiplicity
one.
\end{lemma}
\begin{proof}
The proof follows from  Bertini's Theorem applied to the linear
system defining $\mathcal G(\mathcal F)$, or equivalently, from
Sard's Theorem applied to $\mathcal G(\mathcal F)$. For the details
see \cite{CLS}.
\end{proof}

Notice that the conclusion of Lemma concerning the multiplicities
can be rephrased by saying that $H$ is a regular value of  $\mathcal
G(\mathcal F)$ restricted to its domain of definition.

\subsection{Degrees of the Gauss Map}

For a  rational map  $\phi:\mathbb P^n \dashrightarrow \mathbb P^n$
recall that   $\deg_i(\phi)$ is the cardinality of
$\overline{\phi_{|U}^{-1}(L_i)} \cap \Sigma^{n-i}$,
 where $U \subset \mathbb P^n$ is a
Zariski open set where $\phi$ is regular,   $L_i \subset \mathbb
P^n$ is a generic  linear subspace  of dimension $i$ of the target
and $\Sigma^{n-i}\subset \mathbb P^n$ is generic  linear subspace of
 dimension $n-i$ of the domain.

On the remaining part of this section we will study the numbers
${e_i^k(\mathcal F)}$, for pairs of natural numbers $(k,i)$
satisfying $1\le k \le n$ and    $0\le i \le k-1$,   defined as
\[
e_i^k (\mathcal F) =
                                   \deg_i( \mathcal G(\mathcal F_{| \mathbb
P^k}) )   .
\]
where $\mathbb P^k \subset \mathbb P^n$ is generic.

Notice that $e^n_0(\mathcal F)$ is equal to the topological degree
of $\mathcal G(\mathcal F)$ and Lemma \ref{L:CLS} implies that
$e^2_0(\mathcal F) = \deg(\mathcal F)$. More generally, for every $0
\le i \le n-1$, $e^n_i(\mathcal F)$ coincides with the degree of the
{$(n-i)$-th} polar class of $\mathcal F$  defined in \cite{Mol}
mimicking the corresponding definition  for projective varieties,
cf. for instance \cite{Piene}.

Our main result concerning the numbers $e_i^k(\mathcal F)$ is the
following.

\begin{theorem}\label{T:Gauss}
If  $\mathcal F$ is a foliation of $\mathbb P^n$ and  $(k,i)$ is a
pair of natural numbers satisfying $2 \le k \le n$ and $1 \le i \le
k-1$ then
\[
e_i^k(\mathcal F) = e_0^{k-i}(\mathcal F) + e_0^{k-i+1}(\mathcal F)
.
\]
\end{theorem}

The corollary below follows immediately from Theorem \ref{T:Gauss}.

\begin{corollary}
For natural numbers $s,k,i$ satisfying  $s \ge 1$, $s +2 \le k \le
n$ and $2 \le i \le k-1$ we have that $$ e^k_i (\mathcal F) =
e_{i-s}^{k-s} (\mathcal F) \, . $$
\end{corollary}

Notice that this is as an analogous of the invariance of the polar
classes of hypersurfaces under hyperplane sections --- a particular
case of \cite[Theorem 4.2]{Piene}.

\subsection{Proof of Theorem  \ref{T:Gauss}}
It clearly suffices to consider the case $k=n$. Set $U= \mathbb P^n
\setminus \sing(\mathcal F)$ and $\mathcal G = \mathcal G(\mathcal
F)_{|U}$.

Let $L^i \subset \check{\mathbb P}^n$ be a generic linear subspace
of dimension $i$,  $V^i = \mathcal G^{-1}(L^i) \subset U$ and
$\Sigma^{n-i-1} = \check{ L}^i$, i.e.,
\[
 \displaystyle{\Sigma^{n-i-1}=  \bigcap_{H \in L^i} H}
 \, .
\]
Thanks to Bertini's Theorem we can assume that  $V^i$ is empty or
smooth of dimension $i$. Moreover, thanks to Lemma \ref{L:CLS},  we
can also assume that all the singularities of $\mathcal
F_{|\Sigma^{n-i-1}}$ contained in $U$ have multiplicity one.

\begin{lemma}
If  $\Sigma^{n-i}$ is a generic projective subspace of dimension
$n-i$ ($i \ge 1$) containing $\Sigma^{n-i-1}$  then
\[
  V^i \cap \Sigma^{n-i} = U \cap \left( \sing( \mathcal F_{|\Sigma^{n-i}} )
  \cup \sing( \mathcal F_{|\Sigma^{n-i-1}} ) \right)  \, .
\]
Moreover $\Sigma^{n-i}$ intersects $V^i$ transversally.
\end{lemma}
\begin{proof}
By definition $ V^i = \left\{ p \in U \, \vert \, T_p \mathcal F
\supseteq \Sigma^{n-i-1} \right\rbrace .
$
Clearly   the points $p\in \Sigma^{n-i-1}$ belonging  to $V^i$
coincides with $\sing(\mathcal F_{|\Sigma^{n-i-1}})$. Similarly a
point $p \in \Sigma^{n-i} \setminus \Sigma^{n-i-1}$ belongs to $V^i$
if, and only if, $T_p\mathcal F$ contains the join of $p$ and
$\Sigma^{n-i-1}$. Since $\mathrm{Join}(p, \Sigma^{n-i-1}) =
\Sigma^{n-i}$ the set theoretical description of $V^i \cap
\Sigma^{n-i}$ follows.

It remains to prove the transversality statement. First take a point
$p \in \Sigma^{n-i-1} \cap V_i$. If for every $\Sigma^{n-i}$
containing $\Sigma^{n-i-1}$ the intersection of $V^i$ with
$\Sigma^{n-i}$ is not transverse then $T_p V^i \cap T_p
\Sigma^{n-i-1} \neq 0$. Without loss of generality we can assume
that $\Sigma^{n-i-1} = \{ x_0 = \ldots = x_i = 0 \}$. In this
situation the variety $V^i$ is defined by the projectivization of
$\{a_{i+1} = \ldots = a_{n} = 0 \}$ where $\omega = \sum_{i=0}^n a_i
dx_i$ is a $1$-form defining $\mathcal F$ on $\mathbb C^{n+1}$.

If $v \in T_p V^i$ then an arbitrary lift  $\overline v$ to $\mathbb
C^{n+1}$  satisfies  $da_j (\overline v) =0$ for every $ i+1 \le j
\le n$.  Since $\mathcal F_{|\Sigma^{n-i-1}}$ is defined by the
$1$-form
\[
 \sum_{j=i+1}^n a_j(0, \ldots, 0, x_{i+1}, \ldots, x_n) dx_j
\]
then it follows that $d\mathcal G(\mathcal F_{|\Sigma^{n-i-1}})_p
\cdot (v) =0$ for every $v \in T_p \Sigma^{n-i-1} \cap T_p V^i$. If
this latter intersection has positive dimension then $m (\mathcal
F_{|\Sigma^{n-i-1}},p) > 1$ contrary to our assumptions. Therefore
for a generic $\Sigma^{n-i} \supseteq \Sigma^{n-i-1}$ the
intersection of $V^i$ with $\Sigma^{n-i}$ along $\Sigma^{n-i-1}$ is
transversal.

Let now $p \in \Sigma^{n-i} \setminus \Sigma^{n-i-1}$. If $G\subset
\mathrm{aut}(\mathbb P^n)$ is the subgroup that preserves
$\Sigma^{n-i-1}$ then $\mathbb P^n \setminus \Sigma^{n-i-1}$ is
$G$-homogeneous. It follows from  the transversality of a generic
$G$-translate (cf. \cite{Kleiman}) that a generic
$\Sigma^{n-i}\supseteq \Sigma^{n-i-1}$ intersects $V^i$
transversally along $\Sigma^{n-i}\setminus \Sigma^{n-i-1}$.
\end{proof}

The Theorem will follow from the Lemma once we  show that the
closure of $V^i$ in $\mathbb P^n$ cannot intersect $\Sigma^{n-i}
\cap \sing(\mathcal F)$.

For a generic $\Sigma^{n-i}\supset \Sigma^{n-i-1}$ it is clear that
$  \overline{V^i}\cap(\Sigma^{n-i} \setminus \Sigma^{n-i-1})  \cap
\sing(\mathcal F)= \emptyset$. One has just to take a $\Sigma^{n-i}$
transversal to $V^i$ with the maximal number of isolated
singularities contained in $U$.

Our argument to ensure that $  \overline{V^i}\cap\Sigma^{n-i-1} \cap
\sing(\mathcal F)= \emptyset$ is more subtle. Let
$\overline{\mathcal G}:X \to \mathbb P^n$ be a resolution of the
rational map $\mathcal G(\mathcal F)$, i.e, $\pi:X \to \mathbb P^n$
is a composition of smooth blow-ups and $\overline{\mathcal G}$ is
define through the commutative diagram below.
\[
 \xymatrix{
 X  \ar[d]_{\pi}  \ar@/^0.4cm/[dr]^{\overline{\mathcal G}}   \\
  \mathbb P^n \ar@{-->}[r]^{\mathcal G(\mathcal F)} & \check{\mathbb P}^n}
\]
Let also $\mathcal I \subset \mathbb P^n \times \check{\mathbb P}^n$
be the incidence variety, $ \mathbb G_{i}(\check{\mathbb P}^n)$ be
the Grassmanian of  $i$-dimensional linear subspaces of
$\check{\mathbb P}^n$ and
\[
\mathcal U = \left\{ ( L^i, x , H) \in \mathbb G_{i}(\check{\mathbb
P}^n)\times \mathbb P^n \times \check{\mathbb  P}^n  \, \Big\vert \,
  H \in L^i, x \in \check{L}^i= \bigcap_{H \in L^i} H \right\} .
\]
Notice that $\mathcal U \subset \mathbb G_{i}(\check{\mathbb P}^n)
\times \mathcal I$.

If $E \subset X$ is an exceptional divisor then the set of
$i$-dimensional linear subspaces $L^i \subset \check{\mathbb P}^n$
for which $\overline{\mathcal G^{-1}}(L^i) \cap
\pi^{-1}(\check{L}^i) \cap E \neq \emptyset$ is given by the image
of the morphism $\sigma$  defined below, where the unlabeled
 arrows are the corresponding natural projections.
\[
 \xymatrix{
 E \times_{\mathcal I}\ar@/^0.5cm/@{->}[rr]^{\sigma} \mathcal U \ar[d] \ar[r]  & \mathcal U  \ar[d] \ar[r]  & \mathbb G_i(\check{\mathbb P}^n)
  \\
  E \ar[r]^{\pi \times\overline{ \mathcal G} } & \mathcal I}
\]
Notice that $\mathcal I$ is a $\mathrm{aut}(\mathbb
P^n)$-homogeneous space under the natural action and that the
vertical arrow $\mathcal U \to \mathcal I$ is a
$\mathrm{aut}(\mathbb P^n)$-equivariant morphism. The transversality
of the general translate, cf. \cite{Kleiman}, implies that
\[
 \dim  E \times_{\mathcal I} \mathcal U = \dim E + \dim \mathcal U -
 \dim \mathcal I = \dim \mathbb G_i(\check{\mathbb P}^n) -1 .
\]
It follows that $\sigma$ is not dominant. Repeating the argument for
every exceptional divisor of $\pi$ we obtain an open set contained
in $\mathbb G_i(\check{\mathbb P}^n)$ with the desired property.
This concludes the proof of Theorem \ref{T:Gauss}. \qed

\section{Degrees of the Gauss Map of Logarithmic
Foliations}\label{S:log}

Let $F_1, \ldots, F_k \in \mathbb C[x_0, \ldots, x_n]$ be reduced
homogeneous polynomials. If $\lambda = (\lambda_1, \ldots,
\lambda_k) \in \mathbb C^k$ satisfies
\[
\sum_{i=1}^k \lambda_i \deg(F_i) = 0
\]
then the rational $1$-form on $\mathbb C^{n+1}$
\[
\omega_{\lambda} = \omega(F,\lambda) = \sum_{i=1}^k \lambda_i
\frac{dF_i}{F_i} \,
\]
induces a rational $1$-form on $\mathbb P^n$. Formally it is equal
to the logarithmic derivative of the degree $0$ multi-valued
function $F_1^{\lambda_1} \cdots F_k^{\lambda_k}$. Being
$\omega_{\lambda}$ closed, and in particular integrable, it defines
$\mathcal F_{\lambda}$ a singular holomorphic foliation of $\mathbb
P^n$. The corresponding $1$-form is obtained from $(\prod F_i)
\omega_{\lambda}$ after clearing out the common divisors of its
coefficients. The {\it level sets} of the multi-valued function
$F_1^{\lambda_1} \cdots F_k^{\lambda_k}$ are union of leaves of
$\mathcal F_{\lambda}$.

If the divisor $D$ of $\mathbb P^n$  induced by the zero locus of
the polynomial $\prod F_i$ has at most normal crossing singularities
and all the complex numbers $\lambda_i$ are non zero then the
singular of $\mathcal F_{\lambda}$ has a fairly simple structure,
cf. \cite{CHKS,CSV}, which we recall in the next few lines. It has a
codimension two part corresponding to the singularities of $D$ and a
zero dimensional part away from the support of $D$. To obtain this
description one has just to observe that under the hypothesis the
sheaf $\Omega^1(\log D)$ is a locally free sheaf of rank $n$ and
that the rational $1$-form $\omega_{\lambda}$ has no zeros on a
neighborhood of $|D|$ when interpreted as an element of $\mathrm
H^0(\mathbb P^n ,\Omega^1(\log D))$. Moreover, under these
hypotheses,  the length of the zero dimensional part of the singular
scheme of $\mathcal F_{\lambda}$ is measured by the top Chern class
of $\Omega^1(\log D)$.

\smallskip

In order to extend the above description of $\mathrm{sing}(\mathcal
F_{\lambda})$ to a more general setup let
\[
 \pi : ( X , \pi^* D) \to (\mathbb P^n,D),
\]
be an embedded resolution of $D$, i.e., $\pi$ is a composition of
blow-ups along smooth centers contained in the total transforms of
$D$ and the support of $\pi^*D$ has at most normal crossings
singularities.

Due to the functoriality of logarithmic $1$-forms the pull-back
$\pi^* \omega_{\lambda}$ is a global section of
$\mathrm{H}^0(X,\Omega^1_{X}(\log \pi^*D) )$. To each irreducible
component $E$ of $\pi^*D$ there is a naturally attached complex
number --- the {\bf residue} of $\pi^* \omega_{\lambda}$ --- that
can be defined as
\[
\lambda(E)= \lambda(E,\omega_{\lambda})=\frac{1}{2 \pi i}
\int_{\gamma_i} \pi^*(\omega_\lambda)
\]
where  $\gamma : S^1 \to X \setminus |\pi^* D|$  is a naturally
oriented closed path surrounding the support of $E$. If $E$ is the
strict transform of $V(F_i)$ then, clearly,  $\lambda(E) =
\lambda_i$. More generally one has the following lemma.

\begin{lemma}\label{L:nr}
For every irreducible component  $E \subset X$ of the exceptional
divisor there exists natural numbers $m_1, \ldots, m_k \in \mathbb
N$ such that
\[
\lambda(E) = \sum_{i=1}^k m_i \lambda_i \, .
\]
\end{lemma}
\begin{proof}
Let  $\pi_1 : (\mathcal X_1 , \pi_1^*D) \to (\mathbb P^n,D)$ the
first blow up in the resolution process of  $D$ with center $C_1
\subset D$ and let  $E_1=\pi^*(C_1)$ be the exceptional divisor.

If $D_i = V(F_i)$ and $\widetilde{D}_i$ denotes the strict transform
of $D_i$ then we can write
\[
\pi_1^*D_i=n_i {E_1} +  \widetilde{D}_i \,
\]
where $n_i$ is the natural number measuring the multiplicity of
$V(F_i)$ along $C_1$. Moreover if, over a generic point $p \in
|E_1|$, we take $t$ as a reduced germ of regular function cutting
out $E_1$ then
\[
\pi_1^*(\omega_\lambda) =  \left( \sum_{i} \lambda_i n_i\right)
\frac{dt}{t} + \alpha \, ,
\]
for some closed regular  $1$-form $\alpha$. The proof follows by
induction on the number of blow ups necessary to resolve $D$.
\end{proof}

\begin{definition}
The complex vector $\lambda = (\lambda_1, \ldots, \lambda_k) \in
\mathbb C^k$ is {\bf non resonant} (with respect to $\pi$) if
$\lambda(E) \neq 0$ for every irreducible component $E$ of $\pi^*
D$.
\end{definition}

The arguments of \cite{CHKS,CSV}  yields the following description
of the singular set of $\mathcal F_{\lambda}$ for non resonant
$\lambda$. We reproduce them below thinking on  reader's ease.

\begin{lemma}\label{L:CC}
If $\lambda$ is non resonant then the restriction to the complement
of $D$ of the  singular set of $\mathcal F_{\lambda}$ is
zero-dimensional. Moreover the  length of  the corresponding part of
the singular scheme is $c_n(\Omega_X^1 (\log \pi^* D)) \, . $
\end{lemma}
\begin{proof}
Since $\lambda$ is non resonant the $1$-form $\pi^*
\omega_{\lambda}$, seen as a section of $\Omega^1_X(\log \pi^*D)$,
has no zeros on a neighborhood of $|\pi^* D|$.

Suppose that there exists a positive dimensional component of the
singular set of $\mathcal F_{\lambda}$ not contained in $|\pi^* D|$.
Being the divisor $\pi^* D$ ample  this component has to intersect
the support of $\pi^* D$. This leads to contradiction ensuring that
the singular set of $\mathcal F_{\lambda}$ has no positive
dimensional  components in the complement of $|\pi^*D|$.

The assertion concerning the length of the singular scheme follows
from the fact that $\Omega^1_X(\log \pi^* D)$ is a locally free
sheaf of rank $n$.
\end{proof}

Let  $\Sigma^s \subset \mathbb P^n$ be a generic linear  subspace of
dimension $s$ and denote by $X_s = \pi^{-1}(\Sigma^s)$ and $D_s =
(\pi^* D)_{|X_s}$. It follows from Bertini's Theorem that $X_s$ is
smooth and $D_s$ is a divisor with at most normal crossings.

\begin{proposition}\label{P:cc}
If $\lambda$ is non resonant then
\[
\deg_0 (\mathcal G(\mathcal F_{\lambda}) ) =     c_{n-1}(
\Omega^1_{X_{n-1}}(\log D_{n-1}))
\]
and, for $ 1\le i \le n-1 $
\[
    \deg_{n-i} (\mathcal G(\mathcal F_{\lambda}) ) =     c_{i-1}( \Omega^1_{X_{i-1}}(\log D_{i-1}) )+
                                                         c_{i}( \Omega^1_{X_{i}}(\log
                                                         D_{i}))
                                                         \, .
\]
\end{proposition}
\begin{proof}
If $H \subset \mathbb P^n$ is a generic hyperplane then, according
to Lemma \ref{L:CLS}, $\mathcal G(\mathcal F_{\lambda})^{-1}(H)$
coincides with the isolated singularities of $\mathcal F_{|H}$ that
are not singularities of $\mathcal F$. By choosing $H$ on the
complement of the dual variety of the support of $D$ we can assume
that these isolated singularities are away from the support of $D$.

If $\pi_{n-1} : X_{n-1} \to H$ is the restriction of $\pi:X \to
\mathbb P^n$ to $X_{n-1}$ then $\pi_{n-1}$ is an embedded resolution
of $D_{n-1}$ and, moreover, for every exceptional divisor of $E$
intersecting $\pi^{-1}(H)$ we have that the residue of $\pi_{n-1}^*
({\omega_{\lambda}}_{|H})$ along any irreducible component of  $E
\cap X_{n-1}$  is equal to the residue of $\pi^* \omega_{\lambda}$
along $E$.  Therefore the logarithmic $1$-form
${\omega_{\lambda}}_{|H}$ is non resonant with respect to
$\pi_{n-1}$.

It follows from Lemma \ref{L:CC} that the  sought  number of
isolated singularities is $c_{n-1}( \Omega^1_{X_{n-1}}(\log
D_{n-1}))$. Similar arguments shows  that
\[
e^k_0(\mathcal F_{\lambda}) = c_{k-1} ( \Omega^1_{X_{k-1}}(\log
D_{k-1})).
\]
To conclude one has just  to invoke Theorem \ref{T:Gauss}.

\end{proof}

\section{A Logarithmic Foliation associated to a Polar
Transformation}\label{S:geral}

Consider the multivalued function
\[
 \mathbb F^\lambda= \prod_{i=1}^k F_i^{\lambda_i} : \mathbb P^n \dashrightarrow \mathbb P^1
\]
where  $F_i  \in \mathbb C [x_0,...,x_n]$ is a reduced  homogeneous
polynomial of degree $d_i$ and $\lambda_i \in \mathbb C^*$. The
function $\mathbb F^{\lambda}$ is a homogeneous function of degree
$\deg(\mathbb F^{\lambda})=\sum_{i=1}^k \lambda_i d_i$. If
$\deg(\mathbb F^{\lambda})=0$ then the logarithmic derivative of
$\mathbb F^{\lambda}$ defines a logarithmic foliation of $\mathbb
P^n$ and the associated polar map (see below) coincides with the
Gauss map of this foliation. Although one can in principle use the
results of the previous section to control the degree of this polar
map, in general, is rather difficult to control the singular set of
the corresponding logarithmic foliation without further hypothesis.
Therefore, from now on   we will assume that $\deg (\mathbb
F^{\lambda}) \neq 0$.

Although  $\mathbb F^{\lambda}$ is not an algebraic function it is
still possible to define its polar map as the rational map
\begin{eqnarray*}
 \nabla \mathbb F^\lambda:\mathbb P^n &\dashrightarrow& \mathbb P^n \\
  x &\to& \left[ \frac{\mathbb F^{\lambda}_0 (x)}{\mathbb F^\lambda(x)}: \ldots :\frac{\mathbb F^{\lambda}_n(x)}{\mathbb F^\lambda(x)}\right]
\end{eqnarray*}
where  $\mathbb F^{\lambda}_i$ denotes the partial derivative of
$\mathbb F^\lambda$ with respect to  $x_i$. Notice that when all the
$\lambda_i$'s are natural numbers this rational map coincides with
the polar map defined in the introduction.

Consider the foliation of $\mathbb C^{n+1}$ defined by the
polynomial  $1$-form
\[
\left(\prod_{i=1}^k F_i\right) \frac{d\mathbb F^{\lambda}}{\mathbb
F^{\lambda}} = \left(\prod_{i=1}^k F_i\right) \sum_{i=1}^k \lambda_i
\frac{dF_i}{F_i} \, .
\]
Notice that all the singularities of this foliation are contained in
$V(\prod F_i)$ since Euler's formula implies that
\[
i_R \left(\prod_{i=1}^k F_i\right) \frac{d\mathbb
F^{\lambda}}{\mathbb F^{\lambda}} = \deg(\mathbb F^{\lambda})
\left(\prod F_i\right) \, .
\]

This foliation of $\mathbb C^{n+1}$ naturally extends to a foliation
of $\mathbb P^{n+1}$. If we consider $F_1,\ldots, F_k$ as
polynomials in $\mathbb C[x_0,\ldots,x_n,x_{n+1}]$,
$F_{k+1}=x_{n+1}$ and $\overline \lambda = (\lambda_0,
\ldots,\lambda_n,- \deg(\mathbb F^{\lambda}))$ then it coincides
with the foliation $\mathcal F_{\overline \lambda}$ of the previous
section induced by the $1$-form
\[
\omega_{\overline \lambda} = \frac{d\mathbb F^\lambda}{\mathbb
F^\lambda} - \deg(\mathbb F^{\lambda})  \frac{dx_{n+1}}{x_{n+1}} \,
.
\]

The degrees of the Gauss map of $\mathcal F_{\overline \lambda}$ are
related with those of $\mathbb F^{\lambda}$ by means of the
following Theorem.

\begin{theorem}\label{T:polar}
If the degree of $\mathbb F^{\lambda}$ is not equal to zero then for
$i=0,\ldots, n-1$,
\[
\deg_i(\mathcal G(\mathcal F_{\overline \lambda})) = \deg_i
\left(\nabla \mathbb F^{\lambda}\right) + \deg_{i-1} \left(\nabla
\mathbb F^{\lambda}\right) ,
\]
where we are assuming that $\deg_{-1} \left(\nabla \mathbb
F^{\lambda}\right)=0$.
\end{theorem}
\begin{proof}
If we set $\hat F_j = \prod_{i\neq j,i=1}^k F_i$ then the  Gauss map
of the foliation $\mathcal F_{\overline{\lambda}}$ at the point $[
x_0: \ldots:x_{n+1}]$ can be explicitly written as
\[
   \left[ x_{n+1} \left(\sum_{j=1}^k  \lambda_j \hat F_j \frac{\partial  F_j}{\partial
  x_0}\right)
   :\ldots:   x_{n+1} \left(\sum_{j=1}^k \lambda_j \hat F_j \frac{\partial  F_j}{\partial x_n}\right) : -\deg(\mathbb F^{\lambda}) \left(\prod_{j=1}^k F_j
\right)   \right] \, .
\]
Therefore if $p=[0:\ldots:0:1]$ and $\pi:\mathrm{Bl}_p(\mathbb
P^{n+1}) \to \mathbb P^{n+1}$ denotes the blow-up of $\mathbb
P^{n+1}$ at $p$ then the restriction of $\mathcal G= \mathcal
G({\mathcal F_{\overline \lambda}}) \circ \pi^{-1}$ to the
exceptional divisor $E\cong \mathbb P^n$ can be identified with
$\nabla \mathbb F^\lambda$ as soon as we identify the target of
$\nabla \mathbb F^\lambda$ with  the  set  $\mathbb P^n_p\subset
\check{ \mathbb P}^{n+1}$ of hyperplanes containing $p$ .

Consider the projection   $\rho([x_0: \ldots : x_n : x_{n+1}] ) =
[x_0 : \ldots : x_n]$ with center at $p$ and notice that it lifts to
a morphism $\widetilde \rho : \mathrm{Bl}_p(\mathbb P^{n+1}) \to
\mathbb P^n$. If we write
\[
  \nabla \mathbb F^{\lambda} (x) =  \left[ \sum_{j=1}^k  \lambda_j \hat F_j \frac{\partial  F_j}{\partial
  x_0}
   :\ldots:   \sum_{j=1}^k \lambda_j \hat F_j \frac{\partial  F_j}{\partial x_n}   \right] \,
   ,
\]
then it is clear that the rational maps $\mathcal G$ and $\nabla
\mathbb F^{\lambda}$ fit  in the commutative diagram below.
\[
\xymatrix { \mathrm{Bl}_p(\mathbb P^{n+1}) \ar[d]_{\widetilde{\rho}}
\ar@{-->}[rrr]^{\mathcal G} &&& {\mathbb P}^{n+1}
 \ar@{-->}[d]_{{\rho}}\\
\mathbb P^n  \ar@{-->}[rrr]^{\nabla \mathbb F^{\lambda}} &&&
{\mathbb P}^n}
\]

Let $L^i \subset \check{ \mathbb P}^{n+1}$ be a generic linear
subspace of dimension $i$ and  set   $$W^i =\overline{\mathcal
G({\mathcal
 F_{\overline{\lambda}}})^{-1}(L^i)}, \quad  \widetilde{W^i} = \overline{\mathcal G ^{-1}(L^i)} \text{
 and  }
V^i =\overline{ \left(\nabla
  \mathbb F^{\lambda}\right)^{ -1}(\rho(L^i))}.$$ If $U\subset \mathbb P^n$ is the
complement of the hypersurface cut out by $\prod F_j$ then
\cite[lemma]{Piene} implies that
  $V^i\cap U$ and $\widetilde{W^i}
\cap \widetilde{\rho}^{-1}(U)$  are  dense in $V^i$ and
$\widetilde{W^i}$.

It follows at once from the diagram above that
$\widetilde{\rho}(\widetilde{W^i}) \subset V^i$. A simple
computation shows that  the restriction of $\mathcal G$ to a fiber
of $\widetilde{\rho}$ over $U$ induces an isomorphisms to the
corresponding fiber of $\rho$. Combining this  with the density of
$V^i\cap U$ and    $\widetilde{W^i} \cap \widetilde{\rho}^{-1}(U)$
in $V^i$ and $\widetilde{W^i}$ respectively one promptly concludes
that the $i$-cycle $\widetilde{\rho}_* \widetilde{W^i}$ is equal to
the $i$-cycle $V^i$.

   The $i$-th degree of the Gauss map of $\mathcal F_{\overline \lambda}$ can be expressed as
\[
   \deg_i  \left(\mathcal G({\mathcal F_{\overline{\lambda}}})\right) =  c_1(\mathcal O_{\mathbb
   P^{n+1}} (1) )^i \cdot  W^i \, .
\]
If $\widetilde{W^i} = \overline{\mathcal G ^{-1}(L^i)}$, $H$ denotes
a generic hyperplane containing $p$ and $\widetilde H$ is its strict
transform then, thanks to the projection formula,
\begin{eqnarray*}
 \deg_i  \left(\mathcal G_{\mathcal F_{\overline{\lambda}}}\right)
 &=& c_1(\pi^* \mathcal O_{\mathbb
   P^{n+1}} (1) )^i \cdot \widetilde{ W^i}  \\
   &=& c_1( \widetilde\rho^* \mathcal O_{\mathbb P^n}(1))^i \cdot
   \widetilde{W^i} + \left(\sum_{j=1}^i \binom{i}{j}\widetilde H^{i-j} \cdot E^j
   \right)\cdot \widetilde{W^i} \\
   &=&c_1( \mathcal O_{\mathbb
   P^{n}} (1) )^i  \cdot
   \widetilde{\rho} (\widetilde{W^i}) + \left( \left(\sum_{j=1}^i \binom{i}{j}\widetilde H^{i-j} \cdot
   E^{j-1}
   \right) \cap E \right) \cdot\left( \widetilde{W^i} \cap E \right)\\ &=& c_1( \mathcal O_{\mathbb
   P^{n}} (1) )^i  \cdot
   V^i + c_1(\mathcal O_E(1))^{i-1} \cdot (\widetilde{W^i} \cap
   E) \, .
\end{eqnarray*}

On the one hand $c_1( \mathcal O_{\mathbb
   P^{n}} (1) )^i  \cdot
   V^i$ is clearly equal to $\deg_i(\nabla
   \mathbb F^\lambda)$. On the other hand $c_1(\mathcal O_E(1))^{i-1} \cdot (\widetilde{W^i} \cap
   E) = \deg_{i-1}(\nabla \mathbb F^{\lambda})$  since, for a generic $L^i$, $\widetilde{W^i}
\cap
   E$ is equal  to $\overline{\mathcal
   G_{|E}^{-1}(L^i \cap \mathbb P^n _p)}$ as an $(i-1)$-cycle on
   $E$. The Theorem follows.
\end{proof}

\begin{corollary}\label{C:deg}
If the degree of $\mathbb F^{\lambda}$ is not equal to zero then
\[
\deg_i(\nabla \mathbb F^{\lambda}) =  e^{n+1-i}_0 (\mathcal
F_{\overline \lambda} )\, .
\]
for $i=0,\ldots, n-1$.
\end{corollary}
\begin{proof}
Follows at once when after comparing  Theorem \ref{T:Gauss} with
Theorem \ref{T:polar}.
\end{proof}

\section{The Main Result: Invariance of the Degrees}\label{S:main}

\begin{theorem}\label{T:main}
Let $\lambda= (\lambda_1,\ldots, \lambda_k)$ be an element of
$\mathbb C^k$  such that
$
  \mathfrak H(\lambda_j) >0$ for some $\mathbb R$-linear map  $\mathfrak H:\mathbb C \to
\mathbb R$ and every $j=1,\ldots, k $. Let also  $F_1,\ldots, F_k$
be irreducible and homogeneous polynomials in $\mathbb C^{n+1}$. If
$\mathbb F^\lambda = \prod F_j^{\lambda_j}$ then
\[
\deg_i\left(\nabla \mathbb F^\lambda\right) = \deg_i\left(\nabla
\left( \prod F_j \right) \right)
\]
for every $i=0\ldots, n-1$.
\end{theorem}
\begin{proof} Let $\mathcal F= \mathcal
F_{\overline \lambda}$ be the foliation on $\mathbb P^{n+1}$
associated to $\mathbb F^{\lambda}$. Corollary \ref{C:deg} implies
that  $\deg_i(\nabla \mathbb F^{\lambda})$ is equal to the degree of
the Gauss map of $\mathcal F_{|\mathbb P^{n+1-i}}$  for a generic
$\mathbb P^{n+1-i} \subset \mathbb P^{n+1}$.

If  $D$ is the divisor of $\mathbb P^n$ associated to  $\prod F_j$
then the intersection in $\mathbb P^{n+1}$  of
$V\left(x_{n+1}\left(\prod F_j\right) \right)$ and a generic
$\mathbb P^{n-i}$ is isomorphic to the union of
 the intersection of $|D|$ with a generic
$\mathbb P^{n-i}\subset \mathbb P^n$ and a generic hyperplane $H$ in
$\mathbb P^{n-i}$.

If $\pi:X\to \mathbb P^{n-i}$ is an embedded resolution of $|D| \cap
\mathbb P^{n-i}$ then Bertini's Theorem implies that it is also an
embedded resolution of the union of $|D| \cap \mathbb P^{n-i}$ with
a generic $H$. Therefore in the computation of $\lambda(E)$ for  an
exceptional divisor of $\pi$  the residue along $H$,  $\lambda(H)= -
\deg(\mathbb F^{\lambda})$, plays no role since $H$ and its strict
transforms do not contain any of the blow-up centers. Thus the
hypothesis on $\lambda$ together with Lemma \ref{L:nr} implies that
$\overline \lambda$ is non-resonant with respect to $\pi$. It
follows from Proposition \ref{P:cc} that
\[
 \deg_0 (\mathcal G(\mathcal F_{|\mathbb P^{n+1-i}})) =
 c_{n-i}(\Omega^1_X(\log( D\cap \mathbb P^{n-i} + H))) \, .
\]
Since the  same arguments implies that  the same formula holds true
for the foliation associated to $\mathbb F=\prod F_j$ the Theorem
follows.
\end{proof}

The hypothesis on $\lambda \in \mathbb C^k$  can be of course
weakened. Lemma \ref{L:nr} ensures  that  there exits finitely many
subvarieties of $\mathbb C^k$ defined by linear equations with
coefficients in $\mathbb N$ that have to be avoided. Outside these
linear varieties the degree of $\nabla \mathbb F^{\lambda}$ is
constant.

The example below shows, for resonant $\lambda$ the degree of the
associated polar map will in general decrease with respect to the
non-resonant ones.

\begin{example}
Let $F_1, \ldots, F_k,F_{k+1} \in \mathbb C[x,y,z]$ be linear forms
such that $F_1, \ldots, F_{k} \in \mathbb C[x,y]$ and $F_{k+1}
\notin \mathbb C[x,y]$. If $\lambda=(\lambda_1, \ldots, \lambda_k,
\lambda_{k+1}) \in (\mathbb C^*)^{k+1}$ is such that
\[
\sum_{i=1}^{k} \lambda_i = 0
\]
and $k\ge 2$ then the rational map $\nabla \mathbb F^{\lambda}$ is
homaloidal, i.e, $\deg \left( \nabla \mathbb F^{\lambda} \right) =
1$.
\end{example}
\begin{proof}
If  $F_{k+2}$ is a generic linear form and $\lambda_{k+2} = -
\sum_{j=1}^{k+1} \lambda_j = - \lambda_{k+1}$ then the proof of
Theorem \ref{T:main} shows that the degree of $ \nabla \mathbb
F^{\lambda}$ is equal to the number of singularities  of the
foliation $\mathcal F$  of $\mathbb P^2$ induced by the $1$-form
\[
 \left( \prod_{j=1}^{k+2} F_j \right) \sum_{j=1}^{k+1} \lambda_j \frac
 {dF_j}{F_j}
 \]
outside $V\left(\prod_{j=1}^{k+2} F_j\right)$.

Notice that $\mathcal F$ has degree $k$ and that
\[
\sum_{p \in \sing(\mathcal F)} m(\mathcal F,p)  =
c_2(\Omega^1_{\mathbb P^2}(k+2)) = k^2+ k + 1 \, .
\]
On the curve cut out by $\prod F_j$ with $(2k+1) + 1$ singularities.
One of them at the confluence of $k$ lines and the other $2k+1$ at
the intersection of exactly two lines. The latter singularities have
all  multiplicity one as a simple local computation shows. The
multiplicity of the latter can be computed using Van den Essen
formula \cite{van} and is equal to  $k^2 - k -1$. Summing up all
these multiplicities one obtains $k^2 + k$. Thus $\deg (\nabla
\mathbb F^{\lambda}) = 1$.
\end{proof}

In the example above if
\[
\sum_{i=1}^{k} \lambda_i \neq 0 \quad \text{ and } \quad
\sum_{i=1}^{k+1} \lambda_i \neq 0
\]
then Van den Essen Formula shows that the multiplicity of the
singularity containing  the $k$ lines is   $(k-1)^2$. Thus the
degree of $\nabla \mathbb F^{\lambda}$ is, under these hypotheses,
$k-1$. The first author have shown that all the homaloidal polar
maps associated to a product of lines with complex weights are of
the form above. A proof will appear elsewhere.

An easy consequence of Theorem \ref{T:main} is the Corollary below.
It would be interesting to replace the maximum on the left hand side
of the inequality by a sum. Indeed \cite[Proposition 5]{Dimca} does
it for the topological degree  under stronger hypothesis.

\begin{corollary}\label{C:prod}
Let $F_1,F_2 \in \mathbb C[x_0, \ldots, x_n]$ be two  homogeneous
polynomials. If $F_1$ and $F_2$ are relatively prime then
\[
\deg_i ( \nabla F_1 \cdot F_2 ) \ge \max \{ \deg_i ( \nabla F_1
 ), \deg_i ( \nabla  F_2 ) \}
\]
for $i=0,\ldots, n-1$.
\end{corollary}
\begin{proof}
Let $\mathcal F_1$ be the foliation of $\mathbb P^{n+1}$ associated to $F_1$ and $\mathcal F_{12}$ the one
associated to $F_1 F_2$. They are, respectively,
induced by the rational $1$-forms on  $\mathbb P^{n+1}$
\[
\omega_1=\frac{dF_1}{F_1} - \deg(F_1) \frac{dx_{n+1}}{x_{n+1}} \quad
\text{and} \quad \omega_{12}=\frac{dF_1}{F_1} + \frac{dF_2}{F_2} -
(\deg(F_1) \deg(F_2)) \frac{dx_{n+1}}{x_{n+1}}.
\]

Let  $H \subset \check{\mathbb P}^{n+1}$  be  a generic hyperplane   and
$\iota:H  \to \mathbb P^{n+1}$ be the inclusion. Recall that $\mathcal G(\mathcal F_1)^{-1}(H)$ consists
of $\deg_0(\mathcal G(\mathcal F_1))$ isolated points corresponding to the singularities
of $\iota^* \omega_1$ contained in $H \setminus V(F_1)$. It follows from the proof of Theorem \ref{T:main} that
we can assume that $\iota^*\omega_{12}$ is non resonant (with respect to a certain resolution).

  If $H$, seen as a point of $\check{\mathbb P}^{n+1}$, avoids
 the closure of the image of $V(F_2)$ under $\mathcal G(\mathcal F_1)$ then  singularities of $\iota^* \omega_1$ contained
in the complement of  $V(F_1)$ are also contained  in the complement of $V(F_1F_2)$.
It follows that for $\epsilon>0 $ small enough the $1$-form $\iota^* ( \omega_1 + \epsilon \omega_{12})$ has
at least $\deg_0(\mathcal G(\mathcal F_1))$ singularities contained in the complement of $V(F_1F_2)$.
Since  we can choose $\epsilon$ in such a way that $\iota^*(\omega_1 + \epsilon \omega_{12})$ is non resonant the induced
foliation has Gauss map with the degree as the Gauss map of $\mathcal F_{12}$.

It follows from Theorem \ref{T:main} that  $\deg_0(\nabla F_1 F_2 )
\ge \deg_0(\nabla F_1)$. Arguing exactly in the same way first with
$F_2$ and then with  linear sections of higher codimensions the
Corollary follows.
\end{proof}

The Corollary above {\it essentially} reduces the problem of
classification of homaloidal polynomials to the classification of
irreducible homaloidal polynomials and irreducible polynomials with
vanishing Hessian. Although, one should not be much optimistic about
generalizing Dolgachev's Classification to higher dimensions.
Already in  $\mathbb P^3$ there are  examples of irreducible
homaloidal polynomials of  arbitrarily high degree, cf. \cite{CRS}.

\end{document}